[1994/12/01]
\documentclass[10pt]{article}
\title{Functors Extending the Kauffman Backet}
\author{John Armstrong}
\date{September 2006}

\usepackage{amsmath,amsthm,amssymb}
\usepackage{diagrams}
\usepackage{graphicx}

\input{epsf}
\newtheorem{Theorem}{Theorem}[section]

\theoremstyle{definition}
\newtheorem{Definition}[Theorem]{Definition}

\DeclareMathOperator{\Trace}{Tr}

\DeclareMathOperator{\rank}{rank}

\newcommand{\Rmod}[1]{#1\textrm{-}\mathbf{mod}}

\newcommand{\FrTang}{\mathcal{F}r\mathcal{T}ang}

\newcommand{\TL}[2]{\mathcal{TL}_{#2}(#1)}

\begin{document}

\maketitle

\begin{abstract}
  In this article we extend evaluations of the Kauffman bracket on regular isotopy classes of knots and links to a variety of functors defined on the category $\FrTang$ of framed tangles.  We show that many such functors exist, and that they correspond up to equivalence to bilinear forms on free, finitely-generated modules over commutative rings $R$.
\end{abstract}

\renewcommand{\sectionmark}[1]{}

\section{Introduction}
The Kauffman bracket \cite{MR1858113} is one of the most basic polynomial knot invariants, and one with an extremely simple definition by skein relations.  Crossings in a diagram are ``smoothed'' in different ways, and these smoothings are weighted in such a way to leave the sum over all smoothings invariant under diagram moves encoding regular isotopies of links in space.

\begin{figure}[b]
  \epsfxsize=\textwidth
  \epsfbox{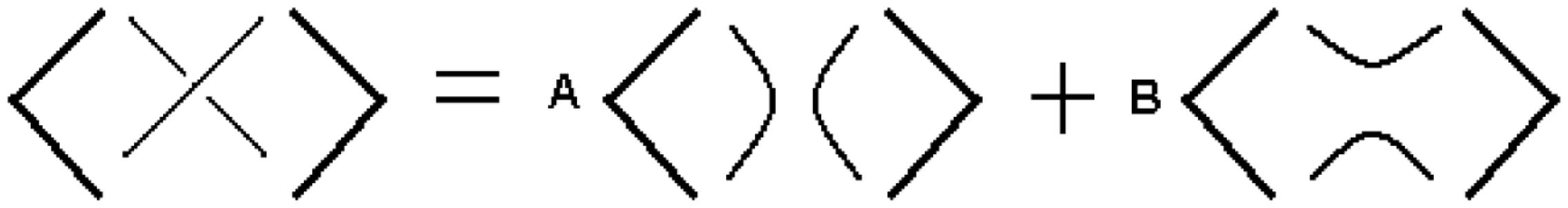}
\caption{The form of the bracket skein relation}\label{fig:BracketSkeinRelation}
\end{figure}

The notion of a skein theory for the bracket is still somewhat fuzzy.  The relation, usually written in the form shown in figure~\ref{fig:BracketSkeinRelation}, is accompanied by words to the effect of ``these diagrams are small parts of larger diagrams which are the same except for around these small parts'', which simply serves to muddle the whole issue.  If we mean to emphasize this area and say that the rest of the link is immaterial, why even consider the rest of the link in the first place?

What we need is to consider the bracket on some sort of ``knots with free ends'', which is exactly what categories of tangles give us.  In particular, since the original bracket is defined on regular isotopy classes of links we will be primarily concerned with the categories $\FrTang$ of framed tangles and $\FrTang(R)$ of $R$-linear combinations of such tangles.  These are monoidal categories which contain framed knots and links as certain morphisms.  We intend to construct functors from these categories which essentially reproduce the bracket when restricted to knots and links.

It is important to distinguish this program of extending knot invariants to functors defined on categories of tangles from the categorification program, and in this case especially from Khovanov's categorification of the bracket \cite{MR1928174}.  Khovanov homology decategorifies to only one functor extending the bracket, and further there is no attempt made in the categorification to respect the monoidal structure which is essential to the simple combinatorial presentations of tangle categories by generators and relations.  In this paper, all functors will be monoidal, and we will identify various families of such functors, in many cases exhausting all possible functorial extensions.

\section{The bracket skein relation on tangles}
\subsection{The categories $\FrTang$ and $\FrTang(R)$}
The category we are most interested in is $\FrTang$: the category of framed, unoriented tangles.  This has a finite presentation \cite{MR1020583}, which we will take as the definition.

\begin{Definition}
  The category $\FrTang$ has as natural numbers as objects with addition as a monoidal structure.  Its morphisms are generated by
  \begin{align*}
    X^+:& 2\rightarrow 2\\
    X^-:& 2\rightarrow 2\\
    \cup:& 0\rightarrow 2\\
    \cap:& 2\rightarrow 0\\
  \end{align*}
  with relations
  \begin{align*}
    T_0:& (\cup\otimes I_1)\circ(I_1\otimes\cap) = I_1 = (I_1\otimes\cup)\circ(\cap\otimes I_1)\\
    T_0':& (I_1\otimes\cup)\circ(X^\pm\otimes I_1) = (\cup\otimes I_1)\circ(I_1\otimes X^\mp)\\
    T_1':& (\cup\otimes\cup)\circ(X^+\otimes X^-)\circ(I_1\otimes\cap\otimes I_1) = \cup\\
    T_2:& X^\pm\circ X^\mp = I_2\\
    T_3:& (X^+\otimes I_1)\circ(I_1\otimes X^+)\circ(X^+\otimes I_1) =\\
    &\qquad\qquad\qquad\qquad= (I_1\otimes X^+)\circ(X^+\otimes I_1)\circ(I_1\otimes X^+)\\
  \end{align*} 
\end{Definition}

A tangle from $m$ to $n$ is interpreted as a tangle diagram in the square with $m$ ends leaving the square on the bottom and $n$ leaving on the top.  The generators $\cup$ and $\cap$ denote a local minimum and maximum between two neighboring strands, respectively, while $X^+$ and $X^-$ denote crossings which would be positive or negative, respectively, if the strands involved were both oriented up the tangle diagram.  Composition is by stacking squares vertically, while the monoidal product is by stacking them horizontally.

Given a commutative ring $R$, we may construct the $R$-linearized category $\FrTang(R)$.  The morphisms from $m$ to $n$ are $R$-linear combinations of $(m,n)$-tangles.  The composition and monoidal product are both extended by $R$-bilinearity.

\subsection{Skein theory and tangles}
Normally a skein relation is presented as a relationship between the values of an invariant at some collection of different links which are ``the same except at one place''.  The various ways that place is to be filled in are drawn as knot diagrams with free ends in the relation.

What we want is to not think of the form of the relation as a sort of shorthand, but to take it seriously.  The diagrams drawn are not to be interpreted as small parts of larger link diagrams, but as diagrams of tangles.  We make the relation between the values of the invariant hold by making it hold at the level of the tangle diagrams themselves.

In terms of tangles, the form of the bracket skein relation is
\begin{equation*}
  \left<X^+\right> = A\left<I_2\right> + B\left<\cap\circ\cup\right>
\end{equation*}
which we will satisfy by demanding
\begin{equation}
  X^+ - A I_2 - B \cap\circ\,\cup = 0
\label{eqn:SkeinRelation1}
\end{equation}

Since this expression involves adding tangles and multiplying them by the weights $A$ and $B$, we should be working within the category $\FrTang(R_{A,B})$, where $R_{A,B}$ is an algebra $R$ over $\Bbb Z[X,Y]$ with $X$ and $Y$ acting as multiplication by $A$ and $B$, respectively.  The bracket (or, rather, a functor extending it) will be taken to preserve whatever $R$-module structure we are working with at the time, and that the target category is $R$-linear as well.

This relation generates an ideal $\mathcal I_{A,B}$ in $\FrTang(R_{A,B})$.  In order that our functors satisfy the relation, we insist that they factor through the quotient category $\FrTang(R_{A,B})/\mathcal I_{A,B}$.  Any such functor is the composition of the canonical projection functor and a functor from $\FrTang(R_{A,B})/\mathcal I_{A,B}$ to the target category.  We must, then, study the structure of this quotient.

\subsection{The category $\FrTang(R_{A,B})/\mathcal I_{A,B}$}
First of all, we can derive a variant of equation~\ref{eqn:SkeinRelation1} by using the other relations in $\FrTang$:
\begin{equation*}
  X^- - B I_2 - A \cap\circ\,\cup = 0
\end{equation*}

Now combining both forms of the skein relation and relation $T_2$, we find that
\begin{align*}
  I_2 &= X^+\circ X^-\\
  &= AB I_2 + A^2 \cap\circ\,\cup + B^2 \cap\circ\,\cup + AB \cap\circ\,\cup\circ\cap\circ\,\cup\\
  &= AB I_2 + (A^2+B^2+AB\bigcirc)\cap\circ\,\cup\\
  (1-AB)I_2 &= (A^2+B^2+AB\bigcirc)\cap\circ\,\cup\\
\end{align*}
where we define $\bigcirc=\cup\circ\cap$.

Now if $1-AB$ is invertible in $R$ (which is always the case if $AB\neq1$ and $R$ is a field) we can use this relation to turn split through any crossing in any way we want, possibly adding a twist to the framing of one strand and multiplying the tangle by some factor in $R$.  For any $m$ and $n$, then, there is only one $R$-parametrized family of tangles for every integer which indicates the total twist of the framing of the tangle.  If $R$ is not a field but $1-AB\neq0$ is not invertible the situation is similar, but a little more combinatorially complex.

In either case, the quotient category is extremely simple unless $AB=1$.  From now on we will assume this is the case and work over $R_A$ where $R$ is an algebra over $\Bbb Z[X,X^{-1}]$ and $X$ acts as multiplication by $A$.  The ideal we are concerned with is now $\mathcal I_A$, generated by $X^+ - A I_2 - A^{-1} \cap\circ\,\cup$.

This implies that
\begin{equation*}
  (A^2+A^{-2}+\bigcirc)\cap\circ\,\cup = 0
\end{equation*}
and so
\begin{equation}
  \bigcirc = -A^2-A^{-2}
\end{equation}
That is, a simple closed loop on its own is equivalent to multiplication by the value $-A^2-A^{-2}$, which we will call $\delta$.

We also check how the skein relation interacts with relations $T_3$, $T_0'$, and $T_1'$.
\begin{align*}
  &(X^+\otimes I_1)\circ(I_1\otimes X^+)\circ(X^+\otimes I_1) =\\
  &A(X^+\otimes I_1)\circ(I_1\otimes X^+)\circ(I_3) + A^{-1}(X^+\otimes I_1)\circ(I_1\otimes X^+)\circ(\cap\circ\cup\otimes I_1) =\\
  &A(X^+\otimes I_1)\circ(I_3)\circ(I_1\otimes X^+) + A^{-1}(I_1\otimes \cup)\circ(\cap\otimes I_1) =\\
  &A (I_3) \circ(X^+\otimes I_1)\circ(I_1\otimes X^+)+ A^{-1}(I_1\otimes\cap\circ\cup)\circ(X^+\otimes I_1)\circ(I_1\otimes X^+) =\\
  &(I_1\otimes X^+)\circ(X^+\otimes I_1)\circ(I_1\otimes X^+)\\
\end{align*}
\begin{align*}
  &(I_1\otimes\cup)\circ(X^+\otimes I_1) =\\
  &A(I_1\otimes\cup)+A^{-1}(\cup\otimes I_1) =\\
  &(\cup\otimes I_1)\circ(I_1\otimes X^\mp)\\
\end{align*}
\begin{align*}
  &(\cup\otimes\cup)\circ(X^+\otimes X^-)\circ(I_1\otimes\cap\otimes I_1)=\\
  &(A^2\bigcirc+\bigcirc\bigcirc+1+A^{-2}\bigcirc)(\cup\otimes\cup)\circ(I_1\otimes\cap\otimes I_1)=\\
  &(-A^4-1+A^4+2+A^{-4}+1-1-A^{-4})\cup=\\
  &\cup\\
\end{align*}

We see that we can reconstruct each of these relations from the skein relation with the help of relations $T_0$ and $T_2$, so they are now superfluous.

In fact, the skein relation can be used to eliminate the crossing generators as well, turning each into a linear combination of tangles involving only cup and cap generators.  We are left with a category of Temperley-Lieb diagrams.

\begin{Definition}
  The $R$-linear category $\TL R\delta$ has as natural numbers as objects with addition as a monoidal structure.  Its morphisms are generated by
  \begin{align*}
    \cup:& 0\rightarrow 2\\
    \cap:& 2\rightarrow 0\\
  \end{align*}
  with relations
  \begin{align*}
    T_0:& (\cup\otimes I_1)\circ(I_1\otimes\cap) = I_1 = (I_1\otimes\cup)\circ(\cap\otimes I_1)\\
    D:& \cup\circ\cap = \delta
  \end{align*} 
\end{Definition}

The composition and monoidal structure in $\TL R\delta$ are analogous to those in $\FrTang(R)$.  In sum, we have proved

\begin{Theorem}
  The category $\FrTang(R_{A})/\mathcal I_{A}$ is equivalent to the category $\TL R{-A^2-A^{-2}}$.
\end{Theorem}

We summarize the various categories and functors involved in the diagram in figure~\ref{fig:SchematicDiagram}.

\begin{figure}[t]
  \begin{diagram}
    \FrTang&&\\
    \dTo(0,2){\iota}&\rdTo(2,2){\left<\_\right>}&\\
    \FrTang(R_A)&\rTo(2,2)&\Rmod R\\
    \dTo(0,4){\pi}&&\uTo(2,2){F}\\
    \FrTang(R_A)/\mathcal I_A&\rTo(2,4){\cong}&\TL R{-A^2-A^{-2}}\\
  \end{diagram}
  \caption{The various categories and functors involved in extending the bracket}\label{fig:SchematicDiagram}
\end{figure}

Here $\iota$ is the canonical inclusion of $\FrTang$ into $\FrTang(R_A)$, $\pi$ is the canonical projection of $\FrTang(R_A)$ onto $\FrTang(R_A)/\mathcal I_A$, and $F$ is a representation of the Temperley-Lieb category, which will be discussed below.

Denoting the functor in the upper-right by $\left<\_\right>$ indicates that when the functor is restricted to regular isotopy classes of knots and links -- morphisms from $0$ to $0$ in $\FrTang$ -- the result is ``essentially'' the bracket of the link.  In fact we will see that the object $0$ must be sent to $R$ as an $R$-module.  The morphism from $R$ to $R$ assigned to any link will be multiplication by the evaluation of the bracket of the link at $A$ in the ring $R$.

\section{Representations of $\TL R\delta$}
The simple presentation of the Temperley-Lieb category and the power of $R$-linear monoidal functors combine to make the representation theory of $\TL R\delta$ much simpler than that of many other categories.  Accordingly, we will be able to state results in a great deal of generality rather than being content to exhibit a handful of examples.

\subsection{The basic structure of a representation}
Firstly, a representation functor $F$ must send the object $1$ in $\TL R\delta$ to some $R$-module $M$.  Since every object in $\TL R\delta$ is some tensor power of $1$, we find that $F$ must send $n$ to $M^{\otimes n}$.  In particular, to preserve the identity of the monoidal structure $F$ must send $0$ to $M^{\otimes 0}\cong R$.

Now the generator $\cap$ must be sent to an $R$-bilinear pairing
\begin{equation*}
F(\cap)=\left<\_,\_\right>:M\times M\rightarrow R
\end{equation*}
while the generator $\cup$ must be sent to a linear ``copairing'' sending $R$ to $M\otimes M$
\begin{equation*}
F(\cup):1\longmapsto \sum x_i\otimes y_i
\end{equation*}
Note that we do not assume either of these have any particular symmetry properties.

The relation $T_0$ tells us that
\begin{equation*}
  \sum x_i\left<y_i,m\right> = m = \sum \left<m,x_i\right>y_i
\end{equation*}
for all $m\in M$.  This indicates that both the $\{x_i\}$ and the $\{y_i\}$ span $M$.  In particular, $M$ must be finitely generated as an $R$-module.

In many cases of interest, $R$ will be a field, in which case $M$ is obviously free.  Even when $R$ is not a field the case of a free $M$ is interesting.  From here on we will take $M$ to be free as well as finitely generated over $R$.  This will also allow us to argue with specific bases of $M$ and matrices for transformations as necessary.

So, let us choose a basis $\left\{e_i\right\}_{i=1}^r$ for $M$, and define the matrices $B=\Bigl(B{i,j}\Bigr)$ and $\hat{B}=\Bigl(\hat{B}{i,j}\Bigr)$ by
\begin{align*}
  F(\cap):&e_i\otimes e_j\longmapsto B{i,j}\\
  F(\cup):&1\longmapsto\sum_{i,j=1}^r \hat{B}{i,j}e_i\otimes e_j\\
\end{align*}

Now relation $T_0$ implies
\begin{equation*}
  \sum_{j=1}^r \hat{B}^{i,j}B_{j,k} = \delta^i_k = \sum_{j=1}^r B_{i,j}\hat{B}^{j,k}
\end{equation*}
These are exactly the equations relating a matrix to its inverse, and we can use Cramer's rule to calculate $\hat{B}$ when $B$ is given.  The copairing is thus determined uniquely by the pairing.  If $R$ is a field then such a copairing will always exist as long as the pairing is nonsingular, while over arbitrary rings the conditions are somewhat stronger.  We call a pairing which has such a uniquely defined copairing ``invertible''.

Finally, we find that the relation $D$ implies
\begin{equation}
  \delta = \sum_{i,j=1}^r \hat{B}^{i,j}B_{i,j} = \Trace\left(B^{-1}B^T\right)
  \label{eqn:TraceofAsymmetryCondition}
\end{equation}
The term in the trace on the right is called the ``asymmetry'' of the bilinear form $B$, and will be of use to us later.

We have proved

\begin{Theorem}
  Given a free, finitely-generated $R$-module $M$ with an invertible pairing $\left<\_,\_\right>:M\otimes M\rightarrow R$ with matrix $B$ whose asymmetry has trace $\delta$ there is a unique representation $F_{M,B}$ of the category $\TL R\delta$ sending $\cap$ to $\left<\_,\_\right>$.  Further, all representations sending $1$ to a free $R$-module are of this form.
\end{Theorem}

If $M$ has rank $n$ as an $R$-module we say that the representation $F_{M,B}$ is of rank $n$.

\subsection{Equivalence of representations}
By definition, two representations $F_{M,B}$ and $F_{N,C}$ are equivalent if and only if there is a natural transformation $\eta:F_{M,B}\rightarrow F_{N,C}$.  This consists of a map $\eta_n:M^{\otimes n}\rightarrow N^{\otimes n}$, and to preserve the monoidal structure we must have $\eta_{m+n} = \eta_m\otimes\eta_n$ for any $m$ and $n$.

For the maps $\eta_n$ to generate a natural transformation it must commute with the images of morphisms under $F_{M,B}$ and $F_{N,C}$.  It is sufficient to enforce this on generators.

For the generator $\cap$ we find that the following diagram must commute
\begin{diagram}
R&\rTo(2,0){id_R}&R\\
\uTo(0,0){F_{M,B}(\cap)}&&\uTo(2,0){F_{N,C}(\cap)}\\
M\otimes M&\rTo(2,2){\eta_1\otimes\eta_1}&N\otimes N\\
\end{diagram}
while for the generator $\cup$ we find that the following diagram must commute
\begin{diagram}
M\otimes M&\rTo(2,2){\eta_1\otimes\eta_1}&N\otimes N\\
\uTo(0,0){F_{M,B}(\cup)}&&\uTo(2,0){F_{N,C}(\cup)}\\
R&\rTo(2,0){id_R}&R\\
\end{diagram}

If we write the pairing on $M$ as $\left<\_,\_\right>_M$ and similarly for that on $N$, and pick bases $\left\{m_i\right\}_{i=1}^{\rank M}$ and $\left\{n_i\right\}_{i=1}^{\rank N}$, we find the two conditions
\begin{align*}
  \left<x,y\right>_M &= \left<\eta_1(x),\eta_1(y)\right>\\
  \sum_{i,j=1}^{\rank N} \hat{C}^{i,j}n_i\otimes n_j &= \sum_{i,j=1}^{\rank M} \hat{B}^{i,j}\eta_1(m_i)\otimes \eta_1(m_j)\\
\end{align*}

From relation $T_0$ we find
\begin{equation*}
  n=\sum_{i,j=1}^{\rank M} \hat{B}^{i,j}\left<n,\eta_1(m_i)\right>\eta_1(m_j)
\end{equation*}
so $\left\{\eta_1(m_i)\right\}$ must span $N$.

Now if the rank of $M$ is strictly greater than that of $N$ there must be some nonzero $m$ with $\eta_1(m)=0$.  Then we have
\begin{equation*}
  \left<x,m\right>_M = \left<\eta_1(x),\eta_1(m)\right> = \left<\eta_1(x),0\right> = 0
\end{equation*}
for all $x\in M$, but this is impossible by the nondegeneracy of the pairing on $M$.  Thus $\rank M = \rank N$.  Since $\eta_1$ sends a basis of $M$ to one of $N$, it is an isomorphism, and the natural transformation $\eta$ is a natural isomorphism.

The naturality of $\eta$ with respect to $\cup$ now follows from naturality with respect to $\cap$, and the only condition on $\eta_1$ is that it preserve the pairing.  We have proved

\begin{Theorem}
  Any natural transformation $\eta:F_{M,B}\rightarrow F_{N,C}$ is a natural isomorphism.  There exists such a natural transformation between two functors $F_{M,B}$ and $F_{N,C}$ if and only if the pairings $\left<\_,\_\right>_M$ and $\left<\_,\_\right>_N$ are equivalent.
\label{thm:EquivalenceofRepresentations}
\end{Theorem}

\section{Examples of bracket-extending functors}
As shown in figure~\ref{fig:SchematicDiagram}, a functor from $\FrTang$ will extend the bracket if we pass first to $\TL R{-A^2-A^{-2}}$ and then apply a representation of that category.  However, this is useless if there are no good representations to use.  Luckily, there are actually quite a few that suit our purposes.

\begin{Theorem}
  Given any ring $R_A$ with an identified unit $A\in R$ there is a representation of rank $n$ of $\TL R{-A^2-A^{-2}}$ for all natural numbers $n\geq2$.
\end{Theorem}
\begin{proof}
  If we can find a rank-$n$ representation for $R=\Bbb Z[A,A^{-1}]$ with $A$ as our unit, then for any other $R_A$ we can evaluate our given representation in the new ring to get a representation of $\TL R{-A^2-A^{-2}}$.  We will build our representations out of only the structure of $\Bbb Z[A,A^{-1}]$.
  
  For $n=2$, It is easily checked that
  \begin{equation*}
    \left(
    \begin{array}{cc}
      1&A+A^{-1}\\
      0&1\\
    \end{array}
    \right)
  \end{equation*}
  satisfies equation~\ref{eqn:TraceofAsymmetryCondition} with $\delta=-A^2-A^{-2}$ and has inverse
  \begin{equation*}
    \left(
    \begin{array}{cc}
      1&-A-A^{-1}\\
      0&1\\
    \end{array}
    \right)
  \end{equation*}
  
  Now, for $n\geq2$, let $B$ be an $n$-by-$n$ matrix satisfying equation~\ref{eqn:TraceofAsymmetryCondition} with $\delta=-A^2-A^{-2}$, and with the 1 as the $(n,n)$ entry of $B^{-1}$.  Consider the matrix
  \begin{equation*}
    \left(
    \begin{array}{cc}
      B&\left(\begin{array}{c}0\\\vdots\\1\\\end{array}\right)\\
      \left(\begin{array}{ccc}0&\cdots&0\\\end{array}\right)&1\\
    \end{array}
    \right)
  \end{equation*}
  with inverse
  \begin{equation*}
    \left(
    \begin{array}{cc}
      B^{-1}&-B^{-1}\left(\begin{array}{c}0\\\vdots\\1\\\end{array}\right)\\
      \left(\begin{array}{ccc}0&\cdots&0\\\end{array}\right)&1\\
    \end{array}
    \right)
  \end{equation*}
  
  We calculate
  \begin{align*}
    &\Trace\left(
    \left(
    \begin{array}{cc}
      B^{-1}&-B^{-1}\left(\begin{array}{c}0\\\vdots\\1\\\end{array}\right)\\
      \left(\begin{array}{ccc}0&\cdots&0\\\end{array}\right)&1\\
    \end{array}
    \right)
    \left(
    \begin{array}{cc}
      B^\top&\left(\begin{array}{c}0\\\vdots\\0\\\end{array}\right)\\
      \left(\begin{array}{ccc}0&\cdots&1\\\end{array}\right)&1\\
    \end{array}
    \right)\right)\\
    &= \Trace\left(
    \begin{array}{cc}
    B^{-1}B^\top+\left(\begin{array}{ccc}0&\cdots&(-B^{-1})^1_n\\\vdots&\ddots&\vdots\\0&\cdots&(-B^{-1})^n_n\\\end{array}\right)&-B^{-1}\left(\begin{array}{c}0\\\vdots\\1\\\end{array}\right)\\
    \left(\begin{array}{ccc}0&\cdots&1\\\end{array}\right)&1\\
    \end{array}
    \right)\\
    &=\Trace(B^{-1}B^\top)-1+1=\Trace(B^{-1}B^\top)=\delta
  \end{align*}
  which gives a $(n+1)$-by-$(n+1)$ solution of equation~\ref{eqn:TraceofAsymmetryCondition} with $\delta=-A^2-A^{-2}$ whose inverse has $(n+1,n+1)$ entry 1.
  
  Therefore by induction, we have for each $n\geq2$ at least one representation of $\TL R{-A^2-A^{-2}}$ of rank $n$.
\end{proof}

If $R$ is a field, more can be said.  We know by theorem~\ref{thm:EquivalenceofRepresentations} that our representations correspond to bilinear forms on free $R$-modules up to equivalence.  Although not as well known as the classification of linear transformations up to similarity, there is actually quite a robust theory of the equivalence of bilinear forms, most clearly described in \cite{MR0347867}.  Most interesting to us is the fact that all the methods to determine whether or not two bilinear forms are equivalent start with the asymmetry of the forms -- the linear transformation whose trace we require to be $\delta$ to represent $\TL R\delta$.

We state the most useful case of these conditions explicitly

\begin{Theorem}
  Let $R$ be an algebraically closed field of characteristic $\neq2$. Then two nondegenerate bilinear forms over $R$ are equivalent if and only if their asymmetries are similar.
\end{Theorem}

That is: given a nondegenerate bilinear form $B$ on a vector space $V$ over $R$ with asymmetry $G$ we can find a basis of $V$ putting $G$ into Jordan normal form.  The Jordan block decomposition breaks $V$ into the direct sum of a number of subspaces on which $B$ is indecomposable.  This decomposition is unique up to the order of the summands, as is the case for Jordan normal form.

To find a representation $F_{V,B}$ of rank $n$ we must find a bilinear form of rank $n$ such that $\Trace(G)=\delta$.  Consider an indecomposable component whose asymmetry is the Jordan block
\begin{equation*}
  \left(
  \begin{array}{cccc}
    \lambda&1&\cdots&0\\
    0&\lambda&\cdots&0\\
    \vdots&\vdots&\ddots&\vdots\\
    0&0&\cdots&\lambda\\
  \end{array}
  \right)
\end{equation*}
This contributes its dimension $d$ to the rank of $B$ along with $d\lambda$ to the trace of the asymmetry of $B$.  If we have a list of the indecomposable forms that actually occur, along with the Jordan blocks of their asymmetries, the problem of finding a representation for a given rank $n$ and a given $\delta$ reduces to a simple combinatorial problem.  However, producing this list is problematic in many cases.

Luckily such a list does exist in characteristic zero \cite{MR2150889}

\begin{Theorem}
  If $B$ is an indecomposable, nondegenerate bilinear form on a vector space $V$ over an algebraically-closed field of characteristic zero, a basis of $V$ can be found so that the matrix of $B$ is one of the following:
  \begin{equation*}
    H_n(\lambda)=
    \left(
    \begin{array}{cc}
      0&I_n\\
      J_n(\lambda)&0\\
    \end{array}
    \right), \lambda\neq(-1)^{n+1}\\
  \end{equation*}
  \begin{equation*}
    \Gamma_n=
    \left(
    \begin{array}{cccccccc}
      0&0&0&0&\cdots&0&0&(-1)^{n-1}\\
      0&0&0&0&\cdots&0&(-1)^{n-2}&(-1)^{n-2}\\
      \vdots&\vdots&\vdots&\vdots&\ddots&\vdots&\vdots&\vdots\\
      0&-1&-1&0&\cdots&0&0&0\\
      1&1&0&0&\cdots&0&0&0\\
    \end{array}
    \right), n\geq1;
  \end{equation*}
where $J_m(\lambda)$ is a Jordan block of rank $m$ with eigenvalue $\lambda$.  These forms are pairwise non-congruent except for the fact that $H_n(\lambda)$ is congruent to $H_n(\lambda^{-1}$.
\end{Theorem}

A simple calculation shows that $H_n(\lambda)$ contributes $2n$ to the rank of a form, and the trace of its asymmetry is $n(\lambda+\lambda^{-1})$, while $\Gamma_n$ contributes $n$ to the rank and $(-1)^{n+1}n$ to the trace of the asymmetry.  To find all the representations up to equivalence of rank $r$ and with a given $\delta$ value, we need to find all combinations of $H_n(\lambda)$ and $\Gamma_n$ where the sum of the ranks is $r$, and the sum of the traces of the asymmetries is $\delta$.

For example, if we work over $\Bbb C$ we can construct a rank $3$ representation either by using three copies of $\Gamma_1$, one copy of $\Gamma_3$, one copy of $\Gamma_2$ and one of $\Gamma_1$, or one copy of $H_1(\lambda)$ and one of $\Gamma_1$.  In the first and second cases we find the trace of the asymmetry is $3$, so those forms can be used to construct representations of $\TL{\Bbb C}3$.  In the third case the trace of the asymmetry is $-1$, so we can use it to construct a representation of $\TL{\Bbb C}{-1}$.  In the final case we find the trace of the asymmetry is $\lambda+\lambda^{-1}+1$, so by choosing
\begin{equation*}
  \lambda = \frac{\delta+1\pm\sqrt{\delta^2+2\delta-3}}{2}
\end{equation*}
we can use this to construct a representation of $\TL{\Bbb C}\delta$ for any value of $\delta$.

\section{Unitary crossings}
The bracket condition turns out in some cases to be surprisingly restrictive.  If we (working over $\Bbb C$) want to find bracket-extending functors which assign a unitary matrix to the crossing -- for example, to give rise to a family of representations of the braid groups that can be implemented on a quantum computer -- it turns out that our choices are extremely few.

\begin{Theorem}
  If $\left<\_\right>$ is a functor extending a bracket evaluation which sends $1$ to the complex vector space $\Bbb C^n$, then $n$ is $2$ and the point of evaluation $A$ is a fourth root of unity.
\end{Theorem}
\begin{proof}
  Let us write the images of the generators as follows:
  \begin{align*}
    \left<\cup\right>=M\qquad&\qquad M:\mathbb{C}\rightarrow\mathbb{C}^n\otimes\mathbb{C}^n\\
    \left<\cap\right>=N\qquad&\qquad N:\mathbb{C}^n\otimes\mathbb{C}^n\rightarrow\mathbb{C}\\
    \left<X^+\right>=R\qquad&\qquad R:\mathbb{C}^n\otimes\mathbb{C}^n\rightarrow\mathbb{C}^n\otimes\mathbb{C}^n\\
    \left<X^-\right>=L\qquad&\qquad L:\mathbb{C}^n\otimes\mathbb{C}^n\rightarrow\mathbb{C}^n\otimes\mathbb{C}^n\\
  \end{align*}

  Since $X^-$ is the inverse of $X^+$ and both are determined in terms of $\cup$ and $\cap$ by the bracket relation, we can write the unitarity condition as
  \begin{align*}
    aN\circ M + a^{-1}I_{\mathbb{C}^n}\otimes I_{\mathbb{C}^n} &=L=R^\dagger=(aI_{\mathbb{C}^n}\otimes I_{\mathbb{C}^n}+a^{-1}N\circ M)^\dagger\\
    &=\bar{a}I_{\mathbb{C}^n}\otimes I_{\mathbb{C}^n} + \bar{a}^{-1}M^\dagger\circ N^\dagger\\
    \Longrightarrow\qquad(a^{-1}-\bar{a})I_{\mathbb{C}^n}\otimes I_{\mathbb{C}^n}&=\bar{a}^{-1}M^\dagger\circ N^\dagger-aN\circ M
  \end{align*}
  
  The left hand side of this equation has rank $n^2$ unless $\bar{A}=A^{-1}$, while the right hand side can have rank at most 2.  Since $n^2$ is at least $4$, we choose $A=e^{i\theta}$ and consider
  \begin{align*}
    M^\dagger\circ N^\dagger &= N\circ M\\
    (M\circ M^\dagger)\circ(N^\dagger\circ N) &= M\circ N\circ M\circ N = (-e^{2i\theta}-e^{-2i\theta})^2=4\cos(2\theta)^2\\
  \end{align*}
  
  If we consider $M$ and $N$ as giving the components of $n$-by-$n$ matrices inverse to each other, we recognize the terms in parens on the left side of the this equation as the squares of the Frobenius norms of $M$ and $N$.  Thus
  \begin{equation*}
    4\geq4\cos(2\theta)^2=\lVert M\rVert^2\lVert N\rVert^2\geq\lVert I_{\mathbb{C}^n}\rVert^2=n^2
  \end{equation*}
  
  So $n$ can only be 2, and $A$ must be $\pm1$ or $\pm i$.
\end{proof}

\section{Conclusions}
We have seen that the invariants on regular isotopy classes of knots and links defined by evaluations of the Kauffman bracket can be seen as the restriction to links of many different functors on the category $\FrTang$ of framed tangles.  The functors that extend the bracket are determined by viewing the skein relation for the bracket not as a functional relation between values of the bracket on links related ``in one place'' in a certain way, but as the generator of an ideal in a category of tangles.

This illustrates the main point of skein theory in the context of tangles: {\em a skein theory for an invariant is a collection of generators for the kernel of the invariant, considered as a functor on a category of tangles}.  We see that the skein-theoretic approach to knot theory is really the study of the ideal theories of categories of tangles, analogously to the way representation theory studies the ideal theories of certain algebras.  From this vantage point, we determine a presentation of a category equivalent to the quotient category through which our invariants must factor, and representations of this category can be determined by specifying the images of its generators subject to constraints given by the relations, just as in the representation theory of groups or algebras.

In the case of the bracket in particular, we see a deep connection between bracket-extending functors and bilinear forms: a form determines a functor, and functors are equivalent if and only if their defining forms are.

\bibliographystyle{amsalpha}
\bibliography{../biblio}
\end{document}